\documentclass[12pt]{article}
\usepackage[all]{xy}

\usepackage{epic}
\usepackage{amsmath}[1996/11/01]
\usepackage{amssymb,amsthm,amsfonts,latexsym,epsfig,here,graphics}
\def\beql#1#2\eeql{\begin{equation}\label{#1}#2\end{equation}}

\advance \textheight by -1 \topmargin
\advance \textheight by -1 \headheight
\advance \textheight by -1 \headsep
\advance \textheight by -2 \footskip
\vsize = \textheight
\hsize = \textwidth
\DeclareMathOperator{\Comp}{Comp}

\DeclareMathOperator{\Aut}{Aut}

\DeclareMathOperator{\SL}{SL}

\newtheorem{theorem}{Theorem}[section]

\newtheorem{kor}[theorem]{Corollary}

\newcommand{\bew}{\noindent\underline{Proof.}\ }

\newtheorem{rem}[theorem]{Remark}
\newtheorem{beisp}[theorem]{Example}
\newtheorem{lemma}[theorem]{Lemma}
\newtheorem{proposition}[theorem]{Proposition}

\newtheorem{defn}[theorem]{Definition}

\renewcommand{\setminus}{-}

\newcommand{\be}{\begin{enumerate}}
\newcommand{\ee}{\end{enumerate}}
\newcommand{\bi}{\begin{itemize}}
\newcommand{\ei}{\end{itemize}}

\newcommand{\ba}{\begin{array}}
\newcommand{\ea}{\end{array}}

\newcommand{\DM}[2]{\mbox{$M _{#1}(#2)^*$}}

\newcommand{\Z}{{\mathbb{Z}}}
\newcommand{\D}{{\mathbb{D}}}
\newcommand{\E}{{\mathbb{E}}}

\newcommand{\F}{{\mathbb{F}}}
\newcommand{\N}{{\mathbb{N}}}

\newcommand{\A}{{\mathbb{A}}}

\newcommand{\C}{{\mathbb{C}}}

\newcommand{\eb}{\phantom{zzz}\hfill{$\square $}\smallskip}

\renewcommand{\em}{\sf}

\begin{document}
\begin{center}
{\Large {\bf  On theta series attached to maximal lattices and their adjoints.}} \\
\vspace{1.5\baselineskip}
{\em Siegfried B\"ocherer}\footnote{Kunzenhof 4B, 79117 Freiburg, Germany,
boecherer@t-online.de} and {\em Gabriele Nebe}
\footnote{Lehrstuhl D f\"ur Mathematik, RWTH Aachen University,
 nebe@math.rwth-aachen.de }
\end{center}

\small
{\sc Abstract.} 
The space spanned by theta series of adjoints of maximal even lattices
of exact level $N$ and determinant $N^2$ 
has the Weierstrass property and hence allows to define extremality
for arbitrary squarefree level $N$. 
We find examples of such dual extremal lattices.
\\
keywords: theta series, modular forms, Weierstrass property, 
dual extremal lattices. 

\normalsize

\section{Introduction}


This paper studies maximal even lattices from the 
geometric, arithmetic and analytic point of view. 
It is interesting to find even lattices $L$ such that the
dual lattice $L^{\#} $ has the highest possible minimum.
The most promising candidates for $L$ are clearly the
maximal even lattices. 

The maximal even lattices $L$ of level $N$ are characterized by the 
arithmetic property that the discriminant group
$L^{\#}/L$ is an anisotropic 
quadratic abelian group of exponent $N$.
If $m:=\dim(L) =2k$ is even, 
then this property can be translated in transformation rules 
of the theta series of $L$ under the Atkin-Lehner involutions
for all prime divisors of $N$ (Theorem \ref{thetamaxpthm}).
If $\det (L) = N^2$ then the theta series of the adjoint lattice 
$\sqrt{N} L^{\#} $ lies in the
space $M_k(N)^{*}$ introduced in \cite{A-B}. 
This space has the Weierstrass property as defined in Definition \ref{Weier}
and hence allows to define extremality. 
The even lattice $L$ is called dual extremal
if the 
 theta series $\theta(\sqrt{N} L^{\#} )$ of the adjoint lattice 
is the extremal modular form in 
$M_k(N)^*$.
The dual extremal lattices of level $N$ 
are the maximal even lattices of level $N$ for which the minimum of
the adjoint lattice is $\geq 2 \dim (M_k(N)^*) $.
Remark \ref{level11} shows that in general this inequality 
may be strict. 
The dimension of $M_k(N)^* $ is calculated in \cite{A-B}.
It is interesting to note that for $k > 2$ the space 
$M_k(N)^*$ is spanned by theta series of adjoint lattices of even maximal
lattices of level $N$, so this space is as small as it can be
to obtain bounds on the minimum with the theory of modular forms. 

The last section of this paper lists some examples of dual extremal 
lattices. The level 2 case is remarkable.
Its connection to the notion of s-extremal (odd) unimodular lattices
in \cite{Gaborit} allows to prove that 
for a dual extremal lattice $L$ of level $2$ and dimension $2k$
the minimum
$\min (L^{\#} ) = \dim (M_k(2)^*) $. 
Also for $k\equiv _{12} 2$ the layers of $L^{\#} $ and of $L$ all 
form spherical 5-designs (Proposition \ref{5des}) and hence both 
lattices are strongly perfect (see \cite{Venkov}) and therefore local maxima of the 
sphere packing density function.
\section{Preliminaries}
\subsection{Modular forms}
For basic facts about modular forms we refer to \cite{Miyake}.
We denote by $M_k(N)$ and $S_k(N)$ the spaces of modular forms and cusp forms
of weight $k$  for the congruence subgroup 
$\Gamma_0(N)=\{\left(\begin{array}{cc} a & b\\
c & d\end{array}\right) \in \SL _2(\Z ) \,\mid c\equiv _N  0\}$.
Throughout the paper, we assume $N$ to be squarefree.
For $\gamma= \left(\begin{array}{cc} a & b\\
c & d\end{array}\right)$ and any function $f$ on the upper 
half plane ${\mathbb H}$ we define the slash operator $\mid_k$ by
$$\left(f\mid_k\gamma\right)(\tau)= det(\gamma)^{\frac{k}{2}} 
(c\tau+d)^{-k}f(\frac{a\tau+b}{c\tau+d})\qquad (\tau\in {\mathbb H}).$$
For primes $p$ we use the Hecke operators $T(p)$ (if $p\nmid N$),
and $U(p)$ (for $p\mid N$) acting on $M_k(N)$ in the usual way. 
We also use the operator $V(p)$ defined by
$$ f\longmapsto \left(f\mid V(p)\right)(\tau):= f(p\cdot \tau).
$$
Occasionally we need a variant $U(p)^0$ of the operator $U(p)$, defined for
functions $f$ on ${\mathbb  H}$ periodic with respect to $p\cdot {\mathbb Z}$:
$$f(\tau)=\sum_n a_n e^{2\pi i \frac{n}{p}\tau}\longmapsto f\mid U^0(p)(\tau)=
\sum_n a_{np} e^{2\pi i n\tau}.$$ 
Let $p$ be a prime with $p\mid N$.
We denote by $\omega^N_p$ any element of $SL(2,{\mathbb Z})$ satisfying

$$ \omega_p^N\equiv \left(\begin{array}{cc} 0 & -1\\
1 & 0\end{array}\right)\bmod p $$
and
$$ \omega_p^N\equiv 1_2 \bmod \frac{N}{p}.$$

For such a matrix $\omega_p^N$ we put

$$W_p^N:= \omega_p^N\cdot \left(\begin{array}{cc} p & 0\\
0 & 1\end{array}\right)$$

and we recall that such a matrix defines an ``Atkin-Lehner involution''
on the space $M_k(N)$.

\subsection{Lattices}          

We mainly consider even lattices $L$ in some positive definite
quadratic space $(V,Q)$.
Here $L$ is called {\em even}, if $Q(L) \subset \Z $.
Then $L$ is automatically contained in its dual lattice
$L^{\#} := \{ x\in V \mid (x,\ell ) \in \Z \mbox{ for all } \ell \in L \}$ 
where $(x,y) := Q(x+y) -Q(x) - Q(y)$ is the associated bilinear form. 
The minimal number $N\in \N$ such that the 
{\em adjoint lattice}
$\sqrt{N} L^{\#} := (L^{\#} , N Q ) $ is again even is 
called the {\em level} of $L$. 
We also define the {\em minimum} 
$\min (L) := \min \{ (\ell, \ell) \mid 0\neq \ell \in L \} $.

For a quadratic space $(V,Q)$ over ${\mathbb Q}$ we define the 
local Witt invariants $s_p(V)$ as in \cite[p.80]{Scharlau}. This normalization
is very convenient for our purposes, in particular we will use the
following lemma from 
 \cite{BFSP}.

\begin{lemma}
Let $L$ be an even lattice of level $N\cdot p$ with 
$p\nmid N$ in the quadratic space $(V,Q)$ Then the following statements 
are equivalent  
\\
i) $s_p(V)=1$\\
ii) $V$ carries (even) lattices of level $N$.\\ 
iii) If $L_p=L_p^{(0)} \perp L_p^{(1)}$ denotes the Jordan splitting of
$L_p = L\otimes \Z _p$, then $L_p^{(1)}$ is an orthogonal sum of 
hyperbolic planes.
\end{lemma}

\section{Lattices maximal at $p$ and their theta series }  \label{thetamaxp}

We assume that $L$ is an even lattice in a positive definite 
quadratic space $(V,Q)$
of dimension $m=2k$ . We denote by $N$ the (exact) level of $L$.
We put $D= det(L)$; then $(-1)^kD$ is a discriminant (i.e. it is 
congruent $1$ or $0$ $\bmod 4$) and we denote by $(-1)^kd$ 
the corresponding fundamental discriminant 
(= a discriminant of a quadratic number 
field or equal to 1). Note that $d$ is odd because $N$ is squarefree.

We consider the theta series
$$\theta(L)(\tau):=\sum_{x\in L} e^{2\pi i Q(x)\cdot \tau} = \sum_{x\in L} q^{ Q(x)}$$
for $\tau\in {\mathbb H}$ and $q=e^{2\pi i \tau}$.
Let $p$ be a prime with $p\mid N$.

We recall the  transformation properties of $\theta(L)$ under $\omega_p^N$:

$$\theta(L)\mid_k\omega_p^N = 
\gamma_p(d_p) s_p(V) D_p^{-\frac{1}{2}} \theta(L^{\sharp, p}) $$
Here $L^{\sharp ,p}= L^{\sharp}\cap {\mathbb Z}[\frac{1}{p}]\cdot L$ is the
lattice dualized only at $p$,
$s_p(V)$ is the Witt invariant 
and
$\gamma_p$ depends only on $d_p\cdot ({\mathbb Q}_p^{\times})^2$, 
more precisely, $\gamma_p(1)=1$ and for 
odd primes $p$, $\delta\in {\mathbb Z}_p^{\times}$

$$\gamma_p(\delta)=1,\qquad \gamma_p(\delta\cdot p)
=(\delta_p,p)_p\cdot (-i)^{\frac{p(p-1)}{2}}$$ 
 
For details see \cite[Lemma 8.2]{opusmagnum}, \cite{BFSP},
or in more classical language, \cite{Kitaoka}, for the explicit 
determination of
$\gamma_p$ see \cite{Funke}. We do not need the more complicated 
$\gamma_2$ here.

\begin{theorem}\label{thetamaxpthm}
 Let $p$ be a prime divisor of $N$ with $p\mid\mid N$.
\begin{eqnarray*}
L_p \quad \mbox{is maximal}& \iff& \\
\theta(L)\mid_k \omega^N_p\mid U^o(p)&=& 
-\gamma_p(d) p^{-1} d_p^{\frac{1}{2}} \,\theta(L)\end{eqnarray*}
\end{theorem}

We remark here that the statement of the theorem is local; actually the 
assumption that $N$ is squarefree is not necessary here.

\bew
"$\Longleftarrow$'':
The transformation properties of theta series imply
$$\theta(L)\mid_k\omega_p^N= \gamma_p(d_p)s_p(V) D_p^{-\frac{1}{2}}  
\theta(L^{\sharp,p})$$
Comparing constant terms on the right sides implies
$$s_p(V)=-1, \qquad 
D_p= p^2\cdot d^{-1}_p.$$
In any case, $(V,Q)$ does not carry a $p$-unimodular lattice
and $D_p=p^2$ or $D_p=p$\\
``$\Longrightarrow $''
Suppose that $L_p$ is maximal, 
in particular, $V_p$ does not carry a lattice, which is unimodular (at $p$),
hence $s_p(V)=-1$.
The local lattice $L_p$ has a decomposition
$$L_p=L_p^{(0)}\perp L_p^{(1)}$$
such that $L_p^{(0)}$ is unimodular and 
the  lattice $\sqrt{p}^{-1} L_p^{(1)}$ 
is anisotropic mod $p$ and of rank 1 or 2.
This implies that any vector in 
$L_p^{\sharp}$ with length in ${\mathbb Z}_p$,
is already in the sublattice 
$L_p$, which implies the
global statement
$$\theta(L^{\sharp,p})\mid U^0(p)=\theta(L).$$
Taking into account that $s_p(V_p)=-1$ and using the transformation
formula from above, we therefore obtain
$$\theta(L)\mid_k\omega_p^N\mid U^0(p)
=-\gamma_p(d) D_p^{-\frac{1}{2}} \theta(L)$$
Moreover, $D_p$ is either $p$ or $p^2$, i.e. $D_p=p^2\cdot d_p^{-1}$.
The assertion follows.
\eb
 
\begin{rem}
 We can more generally consider theta series with 
harmonic polynomials of degree $\nu$,
$$\theta_P(L):=\sum_{x\in L} P(x)e^{2\pi i Q(x)\cdot \tau}.$$
Then we obtain again 
$$\theta_P(L)\mid_{k+\nu} \omega^N_p\mid U^o(p)= 
-\gamma_p(d) p^{-1} d_p^{\frac{1}{2}} \,\theta_P(L)$$
provided that $L_p$ is maximal and $p\mid\mid N$. 
\end{rem}

\begin{rem}
 Theorem \ref{thetamaxpthm}   covers all maximal lattices except those
where the fundamental discriminant $d$ is divisible by $2$ 
(where the level $N$ is
divisible by $4$ and $8$ respectively).
\end{rem}

We will mainly consider lattices which are maximal at all primes $p$.
Concerning the existence we state

\begin{proposition}
Suppose that $N$ is squarefree; then there is an
even maximal lattice of even rank $m=2k$ with $\det(L)=N^2$ if and only if
$m\equiv _8  4 $ and the number of prime divisors of $N$ is odd
or $8\mid m$ and the number of prime divisors of $N$ is even. 
\end{proposition}

\bew
 Let $(V,Q)$ be a quadratic space over ${\mathbb Q}$ possibly carrying
such a lattice.
Then we have for finite primes
$$s_p(V) = -1 \iff p\mid N$$
and 
$$s_{\infty}(V)= \left\{\begin{array}{ccc}
-1 &\mbox{ if } & m\equiv _8  4 \\
1 & \mbox{ if }& 8\mid m
\end{array}\right.$$
By the product formula for the Witt invariant, the number of prime divisors
has to be odd ($m\equiv _8 4 $) or even (if $8\mid m$).
In the other direction we prefer to give an explicit construction:
For $N$ squarefree with an odd number of prime divisors, we choose
a maximal order ${\mathcal O}(N)$ in the quaternion algebra over ${\mathbb Q}$
ramified exactly in the primes dividing $N$. We view it as usual 
as quadratic space (with the norm form).
If $m\equiv _8 4$ we may then take ${\mathcal O}(N)\oplus M$
as an example and for $8\mid m$ we take 
${\mathcal O}(N_1)\oplus {\mathcal O}(N_2)\oplus M$.
Here $M$ is an appropriate even unimodular lattice and 
$N=N_1\cdot N_2$ is a decomposition of $N$ into factors with an 
odd number of prime factors. The maximality of these lattices is then 
easily checked locally. 
\eb

\section{The space $M_k(N)^*$} 

\subsection{Definition and basic properties}
\label{DIMENSIONSFORMEL}
The space of interest for us is (for any squarefree $N>1$ and even weight $k$)
$$M_k(N)^*=\{f\in M_k(N)\,\,\mid\,\, \forall p \mid N :
\,
f\mid W^N_p+p^{1-{k\over 2}}f\mid U(p)=0\} .$$ 
The subspace $S_k(N)^*$ of cuspforms in $M_k(N)^*$ 
was investigated in \cite{A-B}.
We recall some properties from there:\\

\noindent
1) The definition may be rephrased in terms of 
the ``trace''-operator (familiar from the
theory of newforms \cite{Li}):
$$\forall p\mid N: \quad \mbox{trace}^N_{N\over p}(f\mid W_p^N)=0$$
We recall that trace$_{\frac{N}{p}}^N:M_k(N)\longrightarrow M_k(\frac{N}{p})$ 
is defined by $f\longmapsto \sum_{\gamma} f\mid_k\gamma$,
where $\gamma$ runs over $\Gamma_0(N)\backslash \Gamma_0(\frac{N}{p})$;
using explicit representatives for the $\gamma$ we obtain the expression
$\mbox{trace}^N_{\frac{N}{p}}(f)= f + p^{1-\frac{k}{2}} f\mid W^N_p\mid U(p)$
\\
2) When we compare the definition of $S_k(N)^*$ with the 
characterization of newforms in terms of traces, 
we see that $S_k(N)^*$ satisfies half of the conditions
describing newforms, see \cite{Li} for details. In particular, the
 space of newforms of level $N$ is contained in $S_k(N)^*$ and
in fact it is easy to 
see from the theory of newforms that each eigenvalue
system for the collection 
$\{T(p)\in End(S_k(N))\mid p \,\, \mbox{coprime to}\,\, N\}$ 
occurs with multiplicity one in $S_k(N)^*$.
More precisely, $S_k(N)^*$ can be built out of the spaces of 
newforms of level $M\mid N$ as follows:\\
For a normalized Hecke eigenform $f = \sum _{n} a_f(n) q^n $ 
in $S_k(M)^{new}$ we put 
$$f^{(N)}(\tau):= \sum_{d\mid \frac{N}{M}} \mu(d) \frac{d a_f(d)}{\sigma_1(d)}
f(d\cdot \tau)$$
By the same reasoning as in \cite{A-B}, section 2.1, remark 2,
we see that this defines an element of $S_k(N)^*$.
We put
$$S_k(M)^{new,N}:={\mathbb C}\{ f^{(N)}_i\},$$
where $f_i$ runs over the normalized Hecke eigenforms in $S_k(M)^{new}$.  
Then
$$S_k(N)^*= \oplus_{M\mid N} S_k(M)^{new,N}.$$
3) 
We computed the dimension of this space
$$\dim S_k(N)^*={(k-1)N\over 12}-{1\over 2}-{1\over 4}
\left({-1\over (k-1)N}\right)
-{1\over 3}\left({-3\over (k-1)N}\right).$$
4) It is easy to see that $S_k(N)^*$ has codimension one in $M_k(N)^*$,
so there is only one Eisenstein series in this space. 
Actually, we can (at least for $k\geq 4$) compute the Eisenstein series
in $M_k(N)^*$ explicitly from the level one Eisenstein series $E_k$
by the same reasoning as above:

$$E_k^{(N)}:=
\sum_{ d\mid N}\mu(d) \frac{d\sigma_{k-1}(d)}{\sigma_1(d)} E_k(d\cdot\tau).   $$

\subsection{The basis problem for $M_k(N)^*$}
We want to span this space $M_k(N)^*$ by appropriate theta series.
In \cite{Bbasis} we already proved that $S_k(N)^{new}$
is always generated by linear combinations of theta series of quadratic forms 
from any fixed
genus of quadratic forms with (exact) level $N$ and determinant $D$ such that
$p^2\mid D$  and $p^m\nmid D$.
The machinery developed in \cite{Bbasis}, section 8 can also be applied to 
oldforms in $M_k(N)$.

\begin{theorem} \label{thetagen}
 Suppose that the data $m=2k>4, N$ admit the existence of 
a genus ${\mathfrak S}$ of maximal lattices of determinant $N^2$ and rank $m$.
Then
$$M_k(N)^*=\Theta({\mathfrak S}^*),$$
where ${\mathfrak S}^*$ is the genus adjoint to ${\mathfrak S}$ and 
$\Theta({\mathfrak S}^*)$ denotes the ${\mathbb C}$-vector space 
generated by the theta series $\theta(L)$, $L\in {\mathfrak S}^*$.
\end{theorem}

The statement above is false for $m=4$ unless $S_k(N)^*=S_k(N)^{new}$, as 
follows from the work of Eichler \cite{Eichler} and Hijikata-Saito \cite{HS}
on the basis problem.   
Anyway, our proof would not work here (because of convergence 
reasons and because here (and only here) the genus of 
maximal lattices is equal to its adjoint genus).

Before we sketch the proof of this theorem, we recall from Theorem \ref{thetamaxpthm}
 that 
the inclusion $$\Theta({\mathfrak S}^*)\subseteq M_k(N)^*$$ holds.
To simplify the exposition, we only consider the case $N=p$.
We have to study the map

$$\Lambda:\left\{\begin{array}{ccc} 
S_k(p)&\longrightarrow &\Theta({\mathfrak S^*}) \\
g&\longmapsto &\sum_i \frac{1}{m(L_i)} <g,\theta(L_i)>\theta(L_i)
\end{array}\right.
$$
Here $m(L)$ is the number of automorphisms of the lattice $L$ and
the $L_i$ run over representatives of the classes in the genus 
(${\mathfrak S}^*$); the bracket $<,>$ denotes the Petersson product for
modular forms. 
It is a general fact (``pullback formulas'' for Eisenstein series) 
that this map can also be described completely 
in terms of Hecke operators, the explicit form of the contribution of the 
bad place $p$ depends on the genus at hand, see \cite{Bbasis}. 

 The case of newforms of level $p$ was discussed in \cite{Bbasis}.

We just have to add  for a Hecke eigenform $f$ of level one 
an explicit description of the map $\Lambda$ for the
two-dimensional space $$M(f):={\mathbb C}\{f, f\mid V(p)\}.$$ Indeed,
it is of the form

$$\left(\begin{array}{c}
\Lambda(f)\\
\Lambda(f\mid V(p))\end{array}\right) = c\cdot L_2(f,2k-2)\cdot 
{\mathcal A}_p\cdot
\left(\begin{array}{c}
f\\
f\mid V(p) \end{array}\right).$$ 

Here $c$ is an unimportant constant, $L_2(f,s)$ denotes the symmetric square
$L$-function attached to $f$ and ${\mathcal A}_p$ is a certain 
$2\times 2$-matrix (involving the ``Satake parameters'' 
$\alpha_p$ and $\beta_p$ of $f$)
which can be computed from \cite{Bbasis}.
The inclusion $\Theta({\mathfrak S}^*)\subseteq M_k(p)^*$ already implies
that the image of   $M(f)$ under $\Lambda$ is at most
one-dimensional. An inspection of ${\mathcal A}_p$ shows that it is always 
different from the zero matrix (i.e. of rank one), 
in other words, $M(f)$ will always be mapped
onto the one-dimensional space ${\mathbb C}\cdot f^{(p)}\subseteq M_k(p)^*$. 
\begin{rem}
 The case of an arbitrary squarefree number $N$ goes along the 
same line (Kronecker products of such $2\times 2$-matrices have then 
to be considered). A more detailed analysis of these matrices ${\mathcal A}_p$
for arbitrary genera ${\mathfrak S}$ will be given elsewhere \cite{Bbasis2}.
\end{rem}

By the same reasoning (or by applying the Fricke involution 
$\left(\begin{array}{cc} 0 & -1\\
N & 0\end{array}\right)$ to both sides of the theorem) we obtain  

\begin{kor}
 Under the same assumptions as in the theorem we have
$$M_k(N)_*=\Theta({\mathfrak S}),$$
where
\begin{eqnarray*}M_k(N)_*&:=&M_k(N)^*\mid_k \left(\begin{array}{cc} 0 & -1 \\
N & 0\end{array}\right)\\ & = & \{f\in M_k(N)\,\mid \, \forall p\mid N: 
\mbox{trace}^N_{\frac{N}{p}}(f)=0\}\end{eqnarray*}
\end{kor}

{\bf Remark:} Both the theorem and the corollary are remarkable because they 
describe precisely the ``old'' part of $\Theta({\mathfrak S}^*)$ and 
$\Theta({\mathfrak S})$. From the point of view of \cite{BFSP} it may be of 
interest to study the trace of such an oldform: We consider the simplest case,
$N=p$ and $f\in S_k(1)$ is a normalized Hecke eigenform.
Then
$$\mbox{trace}^p_1(f^{(p)})= \mbox{trace}_1^p(f- \frac{p}{p+1}a_f(p)f\mid V(p))= \lambda\cdot f$$
with $\lambda= p+1 - \frac{p}{p+1}a_f(p)^2p^{-k+1}$.
By Ramanujan-Petersson (see \cite{Deligne}) 
$|a_f(p)| \leq 2 p^{(k-1)/2} $ and therefore $\lambda$ cannot be zero. 
On the other hand, $f^{(p)}$ is a linear combination of the $\theta(L)$ 
 with $L\in {\mathfrak S}^*$. The trace of such theta series 
is not understood at all,
see \cite{BFSP}. The situation is completely different
for  $f^{(p)}\mid W_p^p\in S_k(p)_*$: 
this function is in $\Theta({\mathfrak  S})$ and the traces
of the theta series are all zero. This fits well with the fact that 
$tr^p_p(f^{(p)}\mid W_p^p)=0$. 

\subsection{$M(N)_*$ as a module over the ring of modular forms of level one}

The orthogonal sum of a maximal lattice with an even unimodular lattice
is again a maximal lattice. This elementary observation corresponds to
fact that    
 $M(N)_* = \bigoplus _k M_k(N)_*$ is a module over the ring of modular 
forms of level one. The corresponding module structure for $M(N)^*$
is defined by multiplying $f\in M_{k }(N)^*$ with $g\mid_{\ell } \left(\begin{array}{cc} 
0 & -1\\
N & 0\end{array}\right)$ for $g\in M_{\ell }(1)$.

It is clear from the dimension formula for $M_k(N)^*$ that the number of 
generators grows with $N$.

In a few cases we can determine the module structure:
From the point of view of lattices, we may consider
the direct sum $\oplus_{k\equiv_4 2}  M_k(p)_*$ as a module over
the graded ring $\oplus_{4\mid k}  M_k(1)$.
For $p=2$ and $p=3$ the module structure 
is already given in 
Chapter 10 of \cite{cliff}. 
For these two primes $p$, the well known construction A 
establishes an isomorphism between $M_k(p)_*$  and 
the space 
spanned by the Hamming weight 
enumerators of maximal doubly-even self-orthogonal codes $C\leq \F_2^{2k}$ 
respectively maximal self-orthogonal codes $C\leq \F_3^{k}$. 
Note that these weight enumerators span 
the space of relative invariants of the associated Clifford Weil group.
For details we refer to \cite[Chapter 10]{cliff}.
To state the relevant result we need one construction.

\begin{rem}
Let $R\neq \{0\} $ be a root lattice that is an orthogonal summand of the
root sublattice of the 24-dimensional even unimodular lattice $L$.
Then the lattice $M:=\{ \ell \in L \mid (\ell,r) = 0 \mbox{ for all } 
r\in R \} $ is a lattice in dimension $24-\dim(R) $ with
$M^{\#}/M \cong R^{\#}/R $.
Though the isometry class of $M$ does depend on the choice 
of $L$, 
its theta series does not and we will denote it 
 by $\theta (\Comp (R)) := \theta (M)$.
\end{rem}

\begin{proposition}
Let $R:= \bigoplus _{k\equiv _4 0} M_k(1) = \C [\theta (\E _8),\theta (\Lambda _{24})] $ denote the ring spanned
by theta series of even unimodular lattices.
\begin{itemize}
\item[(i)]
For any squarefree $N$ and any $i\in \{ 0,1,2,3 \}$ 
the Fricke-involution is an $R$-module isomorphism between 
$\bigoplus _{k\equiv _4 i} M_k(N)^{*} $ and 
$\bigoplus _{k\equiv _4 i} M_k(N)_{*} $.
\item[(ii)] 
The module $\bigoplus _{k\equiv _4 2} M_k(2)_* $  is the free
$R$-module of rank $2$ with basis $(\theta (\D _4)$, $\theta (\Comp (\D _4)) )$.
\item[(iii)]
The module $\bigoplus _{k\equiv _4 2} M_k(3)_* $  is the free
$R$-module of rank $3$ with basis $(\theta (\A_2\perp \A_2)$,
$\theta (\E _6 \perp \E _6)$, $\theta (\Comp (\A _2\perp A_2)) )$.
\end{itemize}
\end{proposition}

\bew
The first statement is clear, the second one is included in \cite[Theorem 10.7.14]{cliff}
and the last one follows from \cite[Corollary 10.7.7]{cliff}.
\eb

From the point of view of modular forms the full space
$\oplus M_k(p)_*$ deserves attention as well 
as a module over the full graded ring of modular forms of level one.

We consider the cases $p=2, p=3$:\\
{\bf The case $p=2$:}
The dimension formulas show that we will certainly need 
$e_2$, $e_4$ as generators, where $e_2$ is the unique Eisenstein 
series in $M_2(2)_*$ and $e_4=E_4^2\mid W^2_2$ is the unique Eisenstein series in 
$M_4(2)_*$.
The dimension formulas show that
$$\dim M_{k-2}(1) + \dim M_{k-4}(1)= \dim M_k(2)_*$$ 
We can further show that the quotient
$\frac{e_4}{e_2}$ is not a meromorphic modular form of weight $2$ for 
$SL(2,{\mathbb Z})$, therefore
\begin{proposition}
 The space $\oplus_{2\mid k} M_k(2)_*$ is 
a free module over the ring of modular forms of level one with basis
$(e_2, e_4)$. 
\end{proposition}

{\bf The case $p=3$: } Again the dimension formulas show that we need at 
least the generators
$e_2, e_4, h_6$, where $e_2$ and $e_4$ again denote the Eisenstein series
in the space $M_2(3)_*$ and $M_4(3)_*$ and $h_6$ is a nonzero element
in the one-dimensional space $S_6(3)=S_6(3)_*$.
An inspection of the Fourier expansions (in the cusps $\infty$ and $0$ )
shows that a nontrivial relation
$$E\cdot e_2+F\cdot e_4 +H\cdot h_6$$
with level one modular forms $E,F,H$ of weights $k-2, k-4 $ and $k-6$ is
not possible. On the other hand, the dimension formula gives the identity
$$\dim M_{k-2}(1) +\dim M_{k-4}(1) + \dim M_{k-6}(1)= \dim M_k(3)_*,$$
therefore we get
\begin{proposition}
 The space $\oplus_{2\mid k} M_k(3)_*$ is 
 a free module over the ring of modular forms of level one with basis
$(e_2, e_4, h_6 )$. 
\end{proposition}


\section{Extremality}

\subsection{Generalities on analytic extremality}

\begin{defn}\label{Weier}
 A subspace ${\mathcal M}\subseteq M_k(N) $ has 
the 
{\em Weierstrass property} $(\mathcal  W)$ if the projection 
${\mathcal M}\longrightarrow 
{\mathbb C}^r$ to the first $r= dim {\mathcal M} $ coefficients of 
the Fourier expansion
$$f=\sum_{n\geq 0} a_n q^n \longmapsto (a_0,a_1,\dots ,a_{d-1})$$ 
is injective. If this holds, the unique element
$$F=F_{\mathcal M} \in {\mathcal M} $$ with Fourier expansion
$$F=1+ \sum_{n\geq d} a_n q^n$$
is called the {\em extremal modular form} in ${\mathcal M}$.
\end{defn}

If ${\mathcal M}$ contains (say, by definition) 
only modular forms with vanishing
Fourier coefficient $a_0$, the definition of ``Weierstrass property'' has
to be modified in the obvious way. Note that $(\mathcal W)$ holds for
${\mathcal M}$ iff $(\mathcal W)$ holds for the cuspidal subspace 
of $\mathcal M$, provided that the codimension of the 
cuspidal part in ${\mathcal M}$ is one.
\\
The notion ``Weierstrass property'' is motivated by 
the connection of this property with $\infty$ being a 
Weierstrass points on 
the modular curve
$X_0(N)$ if ${\mathcal M}=S_2(N)$, see e.g. \cite{Rohrlich}.
\\
Suppose now that we have a lattice $L$
such that $\theta(L)\in \mathcal M$ for a space $\mathcal M$ with 
property ($\mathcal W$). Then we may call
the lattice $L$ analytically extremal with respect to $\mathcal M$ if
$$\theta(L) =F_{\mathcal M}.$$
In particular, such an analytically extremal lattice satisfies
$$\min(L)\geq 2\cdot\dim ({\mathcal M}).$$
In this generality this definition was introduced in \cite{SchSP}.
\\[0.3cm]

Of course these notions only make sense, if we know interesting classes
of such distinguished subspaces ${\cal M}$.
\begin{beisp}
(1) 
Clearly, for any lattice $L$,  the one-dimensional space 
${\mathcal M}:={\mathbb C}\cdot \theta(L)$ has the property
(${\mathcal W}$) and then $L$ is extremal with respect to this space.\\
(2)
 The full space $M_k(1)$ of modular forms of level 1 has the
Weierstrass property and the well-known Leech lattice is then an
${\cal M}=M_{12}(1)$- extremal lattice.\\
(3)
The spaces of modular forms for the Fricke groups considered by 
Quebbemann \cite{Q1,Q2} 
in his work
on modular lattices.
\end{beisp}

\subsection{Analytic extremality with respect to $M_k(N)^*$}
In general, neither the spaces $S_k(N)$ nor $S_k(N)^{new}$ (or versions of it
appropriately enlarged by some Eisenstein series) have the Weierstrass
property. 
In the case of squarefree level $N$
we showed in \cite{A-B} that the intermediate space
 $S_k(N)^*$ 
(and therefore also $M_k(N)^*$) has the property $({\mathcal W})$,  
therefore there is an extremal modular form 
$$F_{N,k}:=F_{ M_k(N)^*}$$
in this case. 
\begin{defn}
A maximal lattice $L$ of level $N$ and determinant $N^{2}$ in 
dimension $m=2k$ is called {\em dual extremal}, if 
$\theta (\sqrt{N} L^{\#} ) = F_{N,k} $.
\end{defn}

\begin{rem}
(1) 
 Our definition allows to define analytic extremality for 
all squarefree levels. This is in contrast to the situation studied by 
Quebbemann \cite{Q1,Q2}. 
\\
(2)
The additional information $\Theta({\mathfrak S}^*)=M_k(N)^*$ is not 
necessary for the definition of dual extremal lattices, but it 
shows that the space ${\mathcal M}=M_k(N)^*$  the smallest one to 
be chosen for adjoints of maximal lattices.
\end{rem}

\begin{beisp} \label{QUAT}
Let $D$ be a rational
 definite quaternion algebra ramified exactly
at the prime $p$. 
Then any lattice $L$ of level $p$ in the quadratic 
space $(D,n)$, where $n$ is the norm form, is a maximal even lattice.
These lattices $L$ are fractional left-ideals for some maximal order in $D$. 
The non principal $L$ satisfy $\min (L) \geq 4$. 
If the class number (the number of isomorphism classes of left-ideals
for a fixed maximal order in $D$) is two, then 
 $\dim M_2(p)^* = 2$  since  $S_2(p)^*=S_2(p)^{new}$ and 
any non-principal $L$ is dual extremal. 
Note that the definite quaternion algebras over 
${\mathbb Q}$ with class number two are classified
by the work of Kirschmer and Voight \cite{KV}:
 $N=p\in \{11, 17,19\}$ and 
$N=2\cdot 3\cdot 5, N= 2\cdot 3\cdot 7, N= 2\cdot 3\cdot 13 , 
N= 2\cdot 5\cdot 7$.
The condition $S_2(N)^*=S_2(N)^{new} $, which is quite special for the case 
$m=4$, is automatically satisfied if $N=p$, but never in 
the other cases of class number two as can be seen by evaluating the dimension
formula for $S_2(N)^*$.
\end{beisp}

\subsection{A remark about extremal modular forms of level $p$ and weight divisible by $p-1$}
\begin{proposition}\label{ext1}
 Let $p$ be  a prime.
Assume that the weight $k$ is divisible by $p-1$.
Then any modular form $f\in M_k(1)$
with Fourier expansion
$$f\equiv  1+ \sum_{n\geq d} a_nq^n \bmod p \qquad (d=\dim M_k(1))$$
 satisfies
$$f\equiv 1 \bmod p$$
\end{proposition}

\begin{kor}\label{korext}
 Let $p\geq 5$ be a prime. Then 
any extremal modular form $g\in M_k(p)^*$ with 
$k$ divisible by 
$p-1$ satisfies 
  $$g\equiv 1\bmod p$$
\end{kor}

\bew (of Proposition \ref{ext1})
There exists a modular form $\mathcal E$ of weight $k$ with 
${\mathcal E} \equiv 1\bmod p$. For $p\geq 5$ we may take an appropriate
power of the Eisenstein series $E_{p-1}$ of weight $p-1$. For $p=2$ or $p=3$
we can take a suitable monom $E_4^{\alpha}\cdot E_{6}^{\beta}$. 
Therefore we can write $f$ as 
$$f={\mathcal E}+ F$$
with $$F=\sum_{n\geq 1} b_nq^n$$
such that the first $d-1$ coefficients  $b_i$ are congruent zero mod $p$.
For $1\leq i\leq d-1$ we  choose $f_j\in S_k(1)$ with integral Fourier
coefficients $c_{i,n}$ such that for §$1\leq i,j\leq d-1$
$$c_{i,j}=\delta_{i,j}$$
Such cusp forms always exist, see e.g. \cite[Theorem 4.4]{Lang}.
 
Then
$$f={\mathcal E}+\sum a_i f_i + H$$
such that the first $d$ Fourier coefficients of $H$ are zero, hence 
$H$ is identically zero.
The assertion follows.
\eb

To prove the corollary  we note that (by \cite{A-B})
$g$ is equivalent mod $p$ to a modular 
form $G\in M_{k+(p-1)(k-1)}(1)$ provided that $p\geq 5$.  
We apply Proposition \ref{ext1} to this $G$.
  
{\bf Remark:} Using a suitable interpretation of the congruence of 
modular forms, 
it is not necessary in the statements above to assume that 
the Fourier coefficients of
the modular forms are rational.
\\
{\bf Remark:} It would be desirable to include the cases $p=2$ and 
$p=3$ in the corollary.

\section{Examples for dual extremal maximal lattices.}

This section lists some examples of dual extremal maximal lattices
of small level $N$ and small dimension $m$. 
For $N=2$ and $N=3$, one may deduce the classification of all
dual extremal lattices from suitable known classifications of 
unimodular lattices.
For the higher levels $N\geq 5$ we use 
Kneser's neighboring method \cite{Kneser} to list the 
whole genus of maximal lattices together with the 
mass formula to double check the completeness of the result.
The computer calculations where performed with MAGMA.
Gram matrices for the new lattices are available in \cite{database}.
\subsection{$N=2$.} 
Let $L$ be
a maximal 2-elementary lattice  of exact level 2 and 
even dimension $m:=\dim(L) = 2k \equiv _8 4$.
Then 
$L$ is the even sublattice of an odd unimodular lattice $M$
and $L^{\#}  = M \cup v+M$ where $2v\in M$ is a characteristic  
vector of $M$, i.e. $(2v,x) \equiv _2 (x,x) $ for all $x\in M$.
If $\mu = \min (M)$ and $4\sigma $ is the minimal norm of
a characteristic vector in $M$, then $4\sigma \equiv _8 m $ and
$$\min (\sqrt{2} L^{\#} )  = 2 \min ( \mu , \sigma )  .$$
Philippe Gaborit proved in \cite{Gaborit} that for $m\neq 23$
$$\mu +\frac{\sigma}{2} \leq 1+\frac{m}{8} \ \ \star.$$
Lattices achieving this bound are called {\em s-extremal}.
We use $\star $ to show that dual extremal lattices 
$L$ satisfy 
$\min (L^{\#} ) =  \lfloor \frac{k+4}{6} \rfloor $.

\begin{proposition}\label{5des}
Let $L$ be a dual-extremal maximal lattice of level $2$ and
dimension $m=24\ell + 4$. 
Then $L^{\#} $ has minimum $1+2\ell $ and all layers of  $L$ and of 
$L^{\# }$ form spherical 5-designs. In
particular $L$ and   $L^{\# }$ are strongly perfect.
If $M$ is one of the three odd unimodular lattices with even sublattice 
$L$, then $M$ is $s$-extremal of minimum $1+2\ell $.
\end{proposition}

\bew
Let $\mu := \min (M)$ and $\sigma := \min (L^{\#} \setminus M) $.
Since $L$ is dual-extremal $\mu $ and $\sigma $ are both 
$\geq 1 + 2\ell $. 
By the bound in \cite{Gaborit} we obtain 
$\mu + \frac{\sigma }{2} \leq \frac{3}{2} + 3 \ell $ hence
$\mu = \sigma = 1+2\ell $. 
The design property follows from the fact that 
$\dim (\DM{12\ell +2}{2} ) = 
\dim (\DM{12\ell +4}{2} ) = 
\dim (\DM{12\ell +6}{2} ) = 2\ell +1 $.
\eb

Similarly we obtain 

\begin{proposition}
Let $L$ be a dual-extremal maximal lattice of level $2$ and
dimension $m=24\ell - 4$
and let $M$ be one of the three odd unimodular lattices with even sublattice 
$L$. Then $M$ is $s$-extremal of minimum $2\ell $.
The minimum of $L^{\#} \setminus M$ is $2\ell +1$  and
the minimal vectors of 
$L^{\# }$ (which are also those of $L$ and those of $M$)
 form a spherical 3-design, which means that $L^{\#} $, $L$ and $M$ are 
all strongly eutactic.
The lattice $M$ is $s$-extremal.
\end{proposition}

\bew
Let $\mu := \min (M)$ and $\sigma := \min (L^{\#} \setminus M) $.
Since $L$ is dual-extremal $\mu $ and $\sigma $ are both 
$\geq 2\ell $. 
Since $\sigma \equiv_2 \frac{m}{4}$ it is odd
$\sigma \geq 2\ell +1 $.
By the bound $\star $ above we obtain 
$\mu + \frac{\sigma }{2} \leq \frac{1}{2} + 3 \ell $ hence
$\mu = 2\ell,$ $\sigma = 1+2\ell $ and $M$ is $s$-extremal.
\eb

\begin{proposition}
Let $L$ be a dual-extremal maximal lattice of level $2$ and
dimension $m=24\ell +12$. Then
$\min (L^{\#} )  = 2\ell +1$.
\end{proposition}

\bew
Let $M$ be one of the three odd unimodular lattices with even sublattice 
$L$.
Let $\mu := \min(M) $ and $\sigma :=  \min (L^{\#} \setminus M) $.
Since $L$ is dual extremal, $\min (\mu,\sigma ) \geq 2 \ell +1 $. 
By Gaborit's bound $\mu + \frac{\sigma}{2} \leq   3\ell + 2 + \frac{1}{2} $. 
If $\min (\mu ,\sigma ) \geq 2\ell + 2 $, then 
$\mu + \frac{\sigma }{2} \geq 3 \ell + 3$ contradicting the bound above.
\eb

\begin{kor}\label{level2}
A dual extremal lattice $L$ of level $2$ and dimension $2k\equiv _8 4$ 
satisfies $\min (\sqrt{2} L^{\#} ) = 2 \lfloor \frac{k+4}{6} \rfloor $.
\end{kor}

\subsubsection{$m=4$} 
Here the root lattice $\D _4$ is the unique maximal 2-elementary 
lattice and dual extremal.
\subsubsection{$m=12$}
The two root lattices $\D_4\perp \E_8$ and $\D_{12}$ are all 
maximal 2-elementary lattices and both are dual extremal.
\subsubsection{$m=20$}
Let $L$ be a maximal 2-elementary lattice of dimension 20.
Then $L\perp \D _4$ is contained in some even unimodular lattice $U$ of
dimension 24. Since $L$ is maximal it is the orthogonal supplement
$\Comp(\D _4)$
of $\D _4$ in $U$ and $L^{\#}  $ is the projection of $U$ 
to $\D_4^{\perp }$. Since $\min (\sqrt{2} L^{\#} )  \geq 4$, all roots of 
$U$ are either in $\D_ 4$ or perpendicular to this sublattice.
Hence $\D_4$ is an orthogonal summand of the root system of $U$, which is
therefore either $\D _4^6$ or $\D_4 \perp \A _5^4$. 
Both lattices $U$ contain a unique $\Aut(U)$-orbit of such sublattices
$\D _4$ yielding the two dual extremal 2-elementary lattices
of dimension 20.
\subsubsection{$m=28$} 
Let $L$ be a maximal 2-elementary lattice of dimension 28 and 
$M$ be an odd unimodular lattice containing $L$.
If $L$ is dual extremal, then $\min (L^{\#} ) \geq 3$ and hence 
$M$ has minimum $3$.
The 28-dimensional unimodular lattices of minimum 3 are all 
classified in \cite{Bacher}. 
There are 38 isometry classes of such lattices, two of which have
a characteristic vector of norm 4. The other 36 lattices give 
rise to 31 even sublattices $L$ which are all dual extremal. 
By Proposition \ref{5des} 
the $6720$ minimal vectors of $L^{\#} $
as well as all layers of $L$ and $L^{\#} $ 
form  spherical 5-designs and hence $L^{\# }$ is a strongly perfect
lattice (see \cite{Venkov}). 
The next dimension where such a phenomenon occurs is 
$m=52$, where $\min (L^{\#} ) = 5$. 
Then any unimodular sublattice $M$ (with even sublattice $L$) 
is an s-extremal lattice of minimum $5$ in the sense of \cite{Gaborit}. 
Up to now, no such lattice is known.

\subsection{$N=3$.} 

A dual extremal lattice $L$ of dimension $m=2k\equiv _8 4$
satisfies 
$\min (\sqrt{3} L^{\#} ) \geq 2 \frac{k+2}{4} $
\subsubsection{$m=4$.}
Here $\A_2\perp \A_2$ is the unique maximal 3-elementary lattice and
this is dual extremal.
\subsubsection{$m=12$.}
The 3-elementary maximal lattices are $\A_2\perp \A_2\perp \E_8$ and
$\E_6\perp \E_6$, the latter is dual extremal.
\subsubsection{$m=20$.}
Let $L$ be a dual extremal 3-elementary lattice of dimension 20.
Then $L\perp \A_2\perp \A_2$ is contained in an even unimodular
lattice of dimension 24. As for $N=2$ the dual extremality of $L$
implies that the root system of $U$ is $\A_2^{12}$ and there is a
unique such lattice $L$.
\subsubsection{$m=28$.}
Let $L$ be a dual extremal 3-elementary lattice of dimension 28
and let $U$ be an even unimodular lattice of dimension 32
containing $L\perp \A_2\perp \A_2$.
Then $\min (L^{\#} ) \geq 8/3 > 2$ implies that $L$ has no roots and
that the root system of $U$ is $\A_2\perp \A_2$.
By \cite{King} the mass of such lattices $U$ is $>41610$ so there 
are more than $72 \cdot 41610$ such lattices.
Every lattice $L\perp \A_2\perp \A_2$ is contained in 8 unimodular lattices,
so it follows from the discussion below that there are at least 
$9\cdot 41610$ dual extremal lattices. 
The lattice $L^{\#} $ is the projection of $U$ to $(\A_2\perp \A_2)^{\perp }$,
so
$$L^{\#} = \{ x \in (\A_2\perp \A_2)^{\perp } \mid 
\mbox{ there is some } z\in (\A_2\perp \A_2)^{\# }
\mbox{ such that } y:= x+z \in U \} .$$
Here we may assume that $z$ is minimal in its class modulo 
$\A_2\perp \A_2$. Then $(z,z) \in \{ 0, \frac{2}{3} , \frac{4}{3} \} $.
If $x\neq 0$ then $(y,y) \geq 4$ and 
$(x,x) = (y,y) - (z,z) \geq 4 - \frac{4}{3} = \frac{8}{3}$.
This shows that for all these lattices $U$
 the orthogonal $L$  of the root sublattice of $U$ 
is dual extremal.

We list these results and the ones found for level $N=5,7,11$
resp. $N=6,10$  in 
the following tables, with lines  labeled by the level $N$ and
columns labeled by the dimension $m$.
Each entry is the triple $(h,h_{ext},min)$ 
giving the class number $h$ of the genus of maximal lattices,
the number $h_{ext}$ of isometry classes of dual extremal maximal
lattices as well as the minimum $\min (\sqrt{N} L^{\#} )$. 
A ``$\cdot $'' instead of $h$ indicates that we did not compute 
the full genus.
Note that for dimension $m=4$, the classification follows from
Example \ref{QUAT}.

$$
\begin{array}{|c|c|c|c|c|}
\hline
m & 4 & 12 & 20 & 28 \\
\hline
N=2 & (1,1,2) & (2,2,2) & (18 , 2,4) & (\cdot, 31 , 6) \\
\hline
N=3 & (1,1,2) & (2,1,4) & (\cdot , 1, 6) & (\cdot , \geq 9\cdot 41610, 8 ) \\
\hline
N=5 & (1,1,2)  & ( 5,2,4) & (329 , 2 , 8) & \\
\hline
N=7 & (1,1,2) & (12,0,8) & & \\
\hline
N=11 & (3,1,4) & (36,2,10) & (\cdot ,\geq 1,20) & \\
\hline
\end{array}
$$

\begin{rem}\label{level11}
It is interesting to note that for level $N=11$ and
dimension $m=20$, the extremal theta series is $
1 + 132q^{10} + 660q^{12} + 1320q^{13} + 2640q^{14} + \ldots $
so any dual extremal lattice $L$ satisfies 
$\min(\sqrt{11}L^{\#}) = 20 > 2 \dim (M_{10}(11)^*) = 18 $.
So Corollary \ref{level2} does not hold in general for
arbitrary levels.
Note that here the $132$ minimal vectors of $L$ form
a spherical $2$-design.
We constructed such a lattice $L$ as the orthogonal
supplement $L = \Comp(D) \leq \Lambda _{24}$ in the
Leech lattice,
 where $D$ is the dual extremal
lattice of level $11$ and dimension $4$.
\end{rem}

$$
\begin{array}{|c|c|c|c|c|}
\hline
m & 8 & 16 \\
\hline
N=6 & (3,1,4) & (45,2,8) \\
\hline
N=10 & (6,1,6) & (228,7,12) \\
\hline
\end{array}
$$

 \end{document}